\documentclass{article}
\usepackage{graphicx}

\usepackage{graphics,amsmath,amssymb}
\usepackage[margin=1in]{geometry}
\usepackage{amsthm}
\usepackage{amsfonts}
\usepackage{latexsym}
\usepackage{color}
\usepackage{verbatim}
\usepackage{hyperref}
\usepackage[font=small]{caption}
\usepackage{forest}
\usepackage{soul}

\usepackage{url}

\theoremstyle{plain}

\theoremstyle{definition}

\theoremstyle{remark}

\newcommand{\Z}[0]{\mathbb{Z}}

\title{SET! From Groups to Games}
\author{Andrey Boris Khesin\\{\small \textit{University of Oxford}}\\{\small \texttt{andrey.khesin@cs.ox.ac.uk}} \and Tanya Khovanova\\{\small \textit{Massachusetts Institute of Technology}}\\{\small \texttt{tanya@math.mit.edu}}}
\date{}

\begin{document}

\maketitle

\begin{abstract}
    The game of SET is one of the best mathematical games ever. It is no wonder that people have tried to generalize it. We discuss existing generalizations of the game of SET to different groups. We concentrate on two types of generalization: a) where a set consists of cards that multiply to the identity; b) where a set consists of three cards that form an arithmetic progression. We finish with a discussion of some properties of the games that influence how enjoyable they are.
\end{abstract}

\section{Introduction}

Everyone loves games. The game of SET is one of the most mathematical games we know. SET is a great tool to teach mathematics and can be a good introduction to combinatorics, number theory, group theory, geometry, and linear algebra.

SET is played with a deck of cards, where three cards form a \textit{set} if they satisfy some rules. The goal is to quickly find sets among the given cards.

People love this game so much that they have tried to generalize it. At the core of the beautiful properties of the game of SET lies a cyclic group of 3 elements. What will happen if we try a group of 4 or 5 elements instead? What if we try any group?

Some generalizations, such as ProSet and EvenQuads~\cite{LAGames, CragerEtAl, Rose}, are available for purchase. Other generalizations are discussed in the Numberphile video~\cite{N}. Even more generalizations are available for download \cite{tse}.

In this paper, we discuss different generalizations of the game of SET, the mathematics behind these generalizations, and the pros and cons of them. We introduce a few novel generalizations, developed jointly by the first author with Andy Tockman and Della Hendrickson.

The paper is structured as follows. In Section~\ref{sec:SET}, we introduce the rules of the game of SET and establish notation. We can view the SET deck as a group, and a set as three elements that sum to zero. In Section~\ref{sec:m2i}, we discuss generalizations where a set is defined as a deck of cards that multiply to the identity. We also introduce the notion of a torsor. Torsors correspond to decks with no identity element, where all cards are equal. For example, the games ProSet and EvenQuads are based on the same group. The latter game is played on a torsor, while the former game is not. 

We can look at the game of SET differently. We can say that a set consists of three elements of a group that form an arithmetic progression. For the cyclic group of three elements that is underlying the game of SET, this definition is equivalent to the usual one. In Section~\ref{sec:ap}, we discuss generalizations based on this definition.
In Section~\ref{sec:good}, we discuss the advantages and disadvantages of different decks and what features are worth considering when designing your own set variant.

\section{The game of SET}
\label{sec:SET}

\begin{figure}[ht!]
    \centering
    \includegraphics[width=0.6\linewidth]{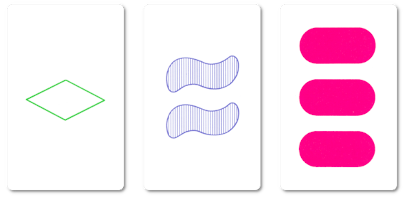}
    \caption{An example of a SET. The underlying group is $C_3^4$.}
\label{fig:set-example}
\end{figure}

The deck for the game of SET consists of 81 cards. Each card has 1, 2, or 3 identical objects. Each object is an oval, a diamond, or a squiggly, colored in red, green, or purple, and the shading can be solid, striped, or empty. Three cards form a set if, for each feature (number, shape, color, shading), all three are either the same or all different. For example, one red striped squiggly, two red striped ovals, and three red striped diamonds form a set. The numbers are all different (1, 2, and 3); the shapes are all different (squiggly, oval, diamond); the color is the same (red), and the shading is the same (striped). Another example of a set with all features different is shown in Figure~\ref{fig:set-example}.

For every feature, we can assign values of 0, 1, or 2, modulo 3. The fact that the values are all the same or all different is equivalent to saying that the values sum up to zero modulo 3. Thus, we can view the cards as points in a 4D vector space over the field $F_3$, namely $F_3^4$. Three cards form a set if and only if the corresponding vectors sum up to the zero vector. Sets have a very useful property: for any two cards in a deck, there exists a third card that completes the given two cards to a set.

Before we generalize the game of SET, we need to establish notation.
Above, we looked at the possible values we assign for a single feature as the elements of the field $F_3$. If we wish to discuss sets in terms of vectors, since a vector space must be over a field, we could use this definition. However, we do not use the multiplicative properties of the field, so we can say that values for one attribute are in $\mathbb{Z}_3$. If we want to forget about the exact values and just use the group structure, we can use the notation $C_3$.

We now discuss possible ways by which we can extend this definition.

\section{Adding or Multiplying to the Identity}
\label{sec:m2i}

Some have tried to generalize the game of SET to a group. We suppose that each card corresponds to an element in a group. Then, we wish to arrange three cards so that the corresponding three elements multiply using the group operation to the group's identity. If the group is non-abelian, then the order of the three cards matters. For the modified game of set, that implies that the players do not just pick out three cards; they have to arrange them in an order so that the corresponding product is the identity. 

However, with this definition, we still retain a useful property. If our card deck contains all possible elements of the group, then for any two cards in a given order, there exists a card that completes them to a set so that the product is the identity.
However, a potential problem is that this card might be one of the two cards we already have. Although it is worth noting that we can potentially rearrange the order of the cards, getting more options.
We could also propose that a set can consist of any number of cards, not just 3, as long as the product of their corresponding group elements is the identity.

The Numberphile video ``The Game of Set (and some variations)'' \cite{N} shows examples of such games. One of the decks they show represents the group $S_3$: permutations of 3 elements. The picture shows a screenshot from the video with a set in a different group, but if you ignore the blue beads, you will get a set in $S_3$. The red dots at the bottom mark odd permutations. This helps see sets faster: each set has to have an even number of red dots. This deck is a toy example of what is possible with only 6 elements in the group.

\begin{figure}[ht!]
    \centering
    \includegraphics[width=0.7\linewidth]{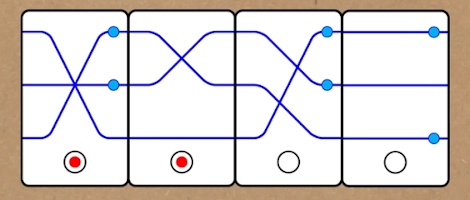}
    \caption{An example set in a game where the four shown group elements multiply to the identity. This is a visual representation of a game played over $C_2\wr S_3$. The red dots (or lack thereof) at the bottom of the cards indicate whether the permutation in $S_3$ is even or odd.}
    \label{fig:three-beaded-wires}
\end{figure}

They also show a similar deck corresponding to group $S_4$: permutations of 4 elements. Thus, the deck size is 24. Despite the size not being very big, the video claims that the game is difficult to play. One more deck represents two permutations of 3 elements, and one must get to the identity on both of them. Thus, the deck size is 36, and the underlying group is $S_3^2$.

These are examples of games over non-abelian groups where the order of the cards forming the set matters. Yet another deck shown is a deck with three crossing lines with beads. The crossing lines correspond to permutations of three elements, and each string either has a bead or not. To get to the identity, one must have an even number of beads on each line. Thus, the deck size is 48 and corresponds to the group that is called the \textit{wreath product} of $C_2$ and $S_3$, where $C_n$ is the cyclic group with $n$ elements. The wreath product is used because a copy of $C_2$, each of the beads, is associated with each of the three lines that the first group, $S_3$, is acting on.

In Figure~\ref{fig:three-beaded-wires}, we see a screenshot from the video, showing four cards forming a set in a game corresponding to this group. The red dot at the bottom signals the parity of a permutation, and there has to be an even number of them, but they can be ignored, as the information they provide is already present on the cards.
However, the meaningful dots are small blue dots on some lines, and there has to be an even number of them on each line.

In practice, one card in each deck will contain the group identity, so this card is typically removed from such games. Leaving it in would allow anyone to add it to any other set they found, getting a free point (if the score is the number of cards for each set). Thus, the actual deck size contains 1 fewer card than the order of the corresponding group.

\begin{figure}[ht!]
    \centering
    \includegraphics[width=0.6\linewidth]{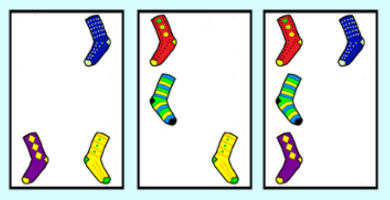}
    \caption{An example set in the game of ProSet, or Socks. Every color of sock appears an even number of times.}
    \label{fig:socks}
\end{figure}

One famous example of such an extension of SET is called ``ProSet'' (short for Projective Set). It is usually marketed under other names, for example, ``Socks''. It is played on the group $C_2^6$ with the identity card removed, resulting in a 63-card deck. Each card contains some non-empty subset of 6 colored dots (or socks). A set is any group of cards where each color appears an even number of times. This is equivalent to adding vectors in $\Z_2^6$ and having them add to the identity: the zero vector. Since the group operation is addition, which is commutative, the cards in the sets found do not need to be in a particular order. A set in the game of Socks is shown in Figure~\ref{fig:socks}.

However, there is a fundamental difference between the classic game of SET and the extensions discussed so far. Unlike in ProSet, no card in the SET deck clearly corresponds to the identity element and thus needs to be removed. It is almost as if the SET group ``forgot'' about its identity element. A group that has forgotten about its identity in this manner is called a \textit{torsor}. While the real definition of a torsor is more sophisticated, this one will do for our purposes.

As discussed, SET sets always contain exactly three cards. These sets have the following property: for each feature, if we assign the numbers 0, 1, and 2 to the three possible values of the feature, we would find that changing this assignment would have the following effect: if the initial set has different values for all of the attributes, those values were assigned the numbers 0, 1, and 2, in some order. A relabeling would merely permute these three values, but all three would still be represented.
If instead the three values for the attributes were all the same, a relabeling would shift each value by the same amount $s$. The sum would thus be shifted by $3s \equiv 0\mod3$.

\begin{figure}[ht!]
    \centering
    \includegraphics[width=0.7\linewidth]{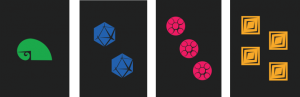}
    \caption{A quad in EvenQuads.}
    \label{fig:EvenQuads}
\end{figure}

There is another game which operates on similar principles, EvenQuads. EvenQuads is played with a deck of size 64. Each EvenQuads card has 3 attributes with 4 values each. Each card has 1 to 4 symbols of one of four colors (red, yellow, green, or blue) and one of four shapes (spiral, polyhedron, circle, or square). A set, which is called a quad in this game, consists of four cards such that for each feature, one of three things is true: 1) the four values for this feature are all the same, 2) the four values are all different, or 3) the values consist of two pairs of identical values (for example, \{circle, circle, square, square\}).
Figure~\ref{fig:EvenQuads} shows a set in this game, where in each feature, all four of the cards are different.

Can we extend the notion of a torsor to ProSet? The answer is yes! EvenQuads is almost equivalent to the game of ProSet \cite{LAGames}. To make ProSet have a torsor-like structure, we just need to add an identity card, only allow sets with four cards, and then map our ProSet cards to those of the EvenQuads deck. The four values of each attribute in EvenQuads can represent all possible values of 2 different socks in ProSet/Socks, which are \{neither sock, one sock, the other sock, both socks\}. With 3 attributes, we can express all possible values of 6 dots/socks.

This representation is more restrictive than the set that we could choose in ProSet, but the advantage is that our deck is now completely free of the choice of a particular identity value. ProSet can be used to play EvenQuads and vice versa (with the exception of the identity card); the ``torsor-ness'' of EvenQuads makes the definition of a set cleaner. Indeed, a set in EvenQuads is formed by any set of four cards that forms a set in ProSet. See more about how different famous card games are equivalent to each other in the paper \textit{Card Games Unveiled: Exploring the Underlying Linear Algebra} \cite{LAGames}.

At this point, it is worth noting that nothing about the games of SET, ProSet, or EvenQuads requires using the exact dimension chosen. The game of SET uses $C_3$ as an underlying group and chooses to use 4 copies of it. Both ProSet and EvenQuads use 6 copies of $C_2$. The choice of dimension, the number of copies, is usually made with practical considerations in mind to make the number of cards in the deck, 81 and 64, respectively, a reasonable quantity. However, to design a new game, we only have to verify that the underlying group is satisfactory, and only then do we need to fix a dimension.

So, how would we generalize further? If each attribute were an element of $\Z_4$ instead of $\Z_2^2$, would that work? Recall that in the game of SET, we require three cards to sum to zero in $\Z_3^4$. In ProSet and EvenQuads, we also require the cards to sum to zero in $\Z_2^6$, where in ProSet, we use any number of cards, and in EvenQuads, we use four cards. Coming back to $\Z_4$, suppose we announce that $x$ cards form a set if their sum is zero, the identity of $C_4$. If we want the cards that are all the same in each feature to form a set, we need $x$ to be divisible by 4. If we choose $x$ to be 4, we will lose the ability to select ``all different'' as a valid combination of attributes since $0+1+2+3 \equiv 2\mod 4$. Any other $x$ that is divisible by 4 is too big to play, although this is only a practical limitation, not a theoretical one. And in any case, we lose symmetry between different values. For example, when the attributes are divided evenly between two values, sometimes they sum to zero, such as in $1+1+3+3\equiv0\mod 4$, and sometimes not, such as in $1+1+2+2\equiv2\mod4$.

Are there larger groups that would work nicely? Let us check $\Z_5$. To allow for all the values to be the same to form a set, we will need sets of 5 cards. The good news is that $0+1+2+3+4 \equiv 0\mod 5$, which means that 5 different values will also satisfy the property of adding to 0. In general, this property will hold for all odd numbers and fail for all even numbers, as the sum of the values from 0 to $n-1$ is $\frac{n(n-1)}{2}$, which is divisible by $n$ if $n$ is odd.

Suppose we define the set as 5 cards that sum to zero modulo 5 across each feature. The big question is, can we define the rule without assigning the exact values to each attribute? Let us look at an example of an acceptable set of values: $\{0, 0, 0, 2, 3\}$. It would not work to always allow sets to contain a 3-1-1 distribution of different values, as that will not always add to 0 mod 5, a requirement we must satisfy if we want the set to be defined on a torsor. However, if we specify a relative order of our elements of the group and allow $\{0, 0, 0, 2, 3\}$ to be a set, then together with it, all shifted sequences correspond to a set; for example, shifting by 1 produces $\{1,1,1,3,4\}$ (note that we still only care about the set of values; their relative order is not relevant).

\begin{figure}[ht!]
\centering
\begin{tikzpicture}
  \node[anchor=south west,inner sep=0] (img) at (0,0)
    {\includegraphics[width=0.8\linewidth]{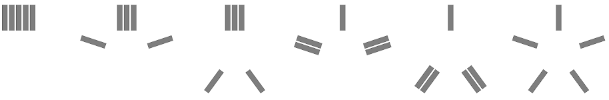}};
  \def\pentsize{0.7cm}
  \def\horizSpacing{0.14225\linewidth}
  \def\vertShift{0.062\linewidth}
  \foreach \i in {0,...,5}{
    \node[draw,regular polygon,regular polygon sides=5,minimum size=\pentsize]
      at($ (img.south west)+(0.024\linewidth+\i*\horizSpacing,\vertShift)$){};
  }
\end{tikzpicture}
\caption{Consider all possible ways to choose five values such that their sum is 0 mod 5. This figure shows all the possible ways that those values can be distributed, up to rotations. As a special feature of the number 5, this is equivalent to choosing all the values that contain an axis of symmetry when placed around a pentagon.}
\label{fig:5Sets}
\end{figure}

In particular, to make a game based on $\Z_5$, we can demand that our sets consist of 5 cards that add to zero, the identity. The fact that each set has 5 cards means that we can use a torsor and do not need to explicitly specify the identity. If we represent our 5 attributes as points on a pentagon or lines pointing in one of 5 directions, we observe a remarkable pattern: the 5 values add to 0 exactly when there is a line of symmetry in the chosen attributes! Figure~\ref{fig:5Sets} shows all possible ways to add 5 values mod 5 and arrive at 0, up to rotations. As an example, if we chose $3+2+1+2+2\equiv0\mod5$, we would see that this corresponds to the second symmetry in the picture.

This is a property unique to the numbers 3 and 5. This is not true for larger sizes of groups; for example, $0+0+0+0+1+2+4=7$ has no symmetries in $\Z_7$. The number 4 represents an interesting case, where all values that add to 0 have symmetries, but not all symmetries are formed by values adding to 0. Specifically, if we represent our 4 attributes as points on a square or lines pointing in one of 4 directions, we observe that when the 4 values add to 0, there is a line of symmetry in the chosen attributes. However, in the cases of the values $\{0, 0, 3, 3\}$ or $\{0, 1, 2, 3\}$, there exists a line of symmetry, but the values do not sum to 0 mod 4.

\begin{figure}[ht!]
    \centering
    
\begin{tikzpicture}
  \node[anchor=south west,inner sep=0] (img) at (0,0) {\includegraphics[width=0.7\linewidth]{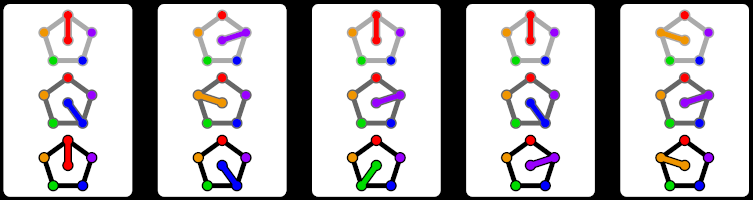}};
  \begin{scope}[shift={(img.south west)},x={(img.south east)},y={(img.north west)}]
    \fill[white,draw=none] (0.183,-0.01) rectangle (0.203,1.01);
    \fill[white,draw=none] (0.3875,-0.01) rectangle (0.4075,1.01);
    \fill[white,draw=none] (0.592,-0.01) rectangle (0.612,1.01);
    \fill[white,draw=none] (0.7965,-0.01) rectangle (0.8165,1.01);
  \end{scope}
\end{tikzpicture}
    \caption{An example of a set in C53T. The first dimension (light gray pentagons) has a vertical axis of symmetry, going through the red vertex. The second dimension (dark gray pentagons) has an axis of symmetry going through the orange vertex. The last dimension (black pentagons) has 5 axes of symmetry. The values across the 3 dimensions can be interpreted as the sums $0+1+0+0+4\equiv2+4+1+2+1\equiv0+2+3+1+4\equiv0\mod5$, although the choice of origin is arbitrary.}
    \label{fig:C53TSet}
\end{figure}
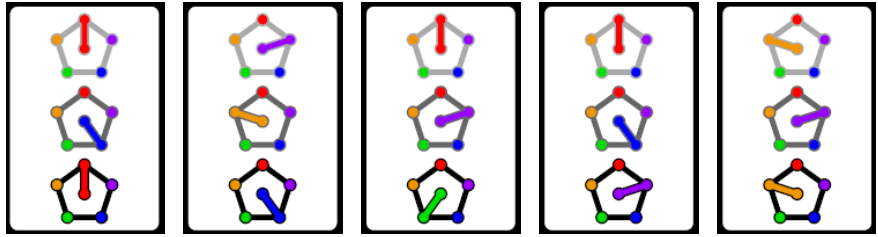

So we have shown that a torsor of $C_5$, the cyclic group of 5 elements, might make for a convenient group on which to play a variant of SET. A reasonable choice of dimension might be 3, making a deck of size 125. This game would then be played on a $C_5^3$-torsor. This inspired the name of this game, which has been shortened to C53T, pronounced ``C-Set''. To illustrate a card, we use pentagons with an indicated vertex to denote an element of $C_5$. We need 3 pentagons, one for each dimension. To help players work around the fact that cards often get turned around when being laid out, each direction on each pentagon is assigned a color to help players tell them apart better. Additionally, the three pentagons are shaded differently to clarify which pentagon corresponds to each dimension. In Figure~\ref{fig:C53TSet}, we see a C53T set, where the symmetric axes in each of the three dimensions go through the red, orange, and ``any'' axes from top to bottom.

This and many other games can be found on the tsetse website \cite{tse}, maintained by Andy Tockman. The website has a description of each game and also allows for solo play.

Finding sets in C53T is extremely hard, so this game is not recommended for casual players, but it does serve to illustrate what really makes a set a set. In everything discussed so far, we have either used a group with a clear identity (such as in the Numberphile video) or a commutative operation (such as in C53T). A natural question that follows is where a non-abelian group torsor can be used to play a variant of SET? As we show in the next section, the answer turns out to be yes.

\section{Sets Forming an Arithmetic Progression}
\label{sec:ap}

Let us return to the original game of SET. In Section \ref{sec:m2i}, we assigned each value one of the numbers 0, 1, or 2, and got the following definition of a set. Three cards form a set if and only if the values for each feature sum to zero modulo 3. Thus, we can see our cards as vectors in the space $\Z_3^4$. Three vectors form a set if they sum to 0.

The generalization we described in Section \ref{sec:m2i} is as follows: We pick a group and define a set as a few cards that might need to be in a specific order that multiply to the group's identity.
For abelian groups, the order does not matter.

However, there is a different way to generalize sets to groups. Three cards that form a set in a classical game of SET, taken in any order, form an arithmetic progression. In other words, if $a$, $b$, and $c$ form a set, then the vectors $b - a$ and $c - b$ are the same. We can check this by considering $a+b+c = a - 2b + c = (a-b) - (b-c)$, which is 0 precisely when $a-b$ and $b-c$ are equal. Here we have used that all vectors are elements of $\Z_3^4$.

Thus, we can generalize the game of SET differently. Suppose our cards are vectors in some space. We say that three of them, $a$, $b$, and $c$, form a set if and only if $b - a = c - b$. Now, the order of the cards in the set becomes important, similar to our previous generalization. In fact, we do not need to use commutative groups like vector spaces. For any group, our condition is equivalent to $ba^{-1} = cb^{-1}$. Thus, for any $a$ and $b$, we have $c = ba^{-1}b$.

Interestingly, this definition does not care about which card is the identity card, meaning the card deck is a torsor. Let us consider a few examples of set variants under this definition.

One potentially desirable feature is to have no group elements of order 2. If this condition is met, there will be no cases when two cards $a$ and $b$ cannot complete to a set because $c=a$. As the game SET already uses the group $C_3^4$, we can consider games based on a power of $C_5$. In fact, we can reuse the cards from C53T for this!

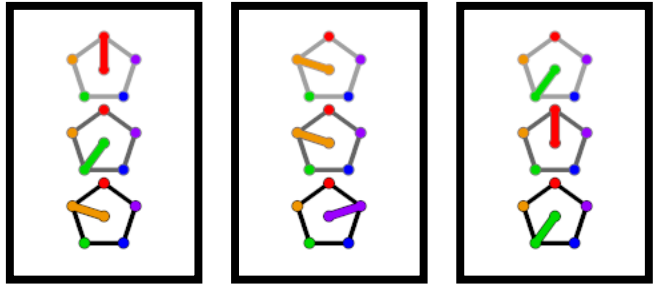
\begin{figure}[ht!]
    \centering
\begin{tikzpicture}[baseline=(img1.base)]
  \node[inner sep=0] (img1) {\includegraphics[scale=0.5,angle=180]{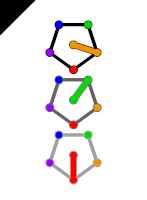}};
  \fill[white] ($(img1.south east)+(-0.8cm,-0.1cm)$) rectangle ++(0.9cm,0.9cm);
  \draw[black,line width=3pt]
    ($(img1.south west)+(0.05cm,-0.05cm)$) rectangle ($(img1.north east)+(-0.05cm,0.05cm)$);
\end{tikzpicture}
\hspace{0.05cm}
\begin{tikzpicture}[baseline=(img2.base)]
  \node[inner sep=0] (img2) {\includegraphics[scale=0.5,angle=180]{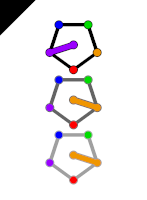}};
  \fill[white] ($(img2.south east)+(-0.8cm,-0.1cm)$) rectangle ++(0.9cm,0.9cm);
  \draw[black,line width=3pt]
    ($(img2.south west)+(0.05cm,-0.05cm)$) rectangle ($(img2.north east)+(-0.05cm,0.05cm)$);
\end{tikzpicture}
\hspace{0.05cm}
\begin{tikzpicture}[baseline=(img3.base)]
  \node[inner sep=0] (img3) {\includegraphics[scale=0.5,angle=180]{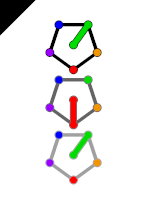}};
  \fill[white] ($(img3.south east)+(-0.8cm,-0.1cm)$) rectangle ++(0.9cm,0.9cm);
  \draw[black,line width=3pt]
    ($(img3.south west)+(0.05cm,-0.05cm)$) rectangle ($(img3.north east)+(-0.05cm,0.05cm)$);
\end{tikzpicture}

    \caption{A set forming an arithmetic progression using C53T cards. The three cards can be interpreted as $(0, 3, 4)$, $(4, 4, 1)$, and $(3, 0, 3)$, although the choice of origin is arbitrary. Both differences are $(4, 1, 2)$, when taken mod 5.}
    \label{fig:linear-C53T-example}
\end{figure}

When playing this set variant with $C_5^3$, three cards form a set if and only if the directions of the first and the third card are symmetric with respect to the direction of the second card. Specifically, we need $c-b=b-a$. An example of such a set is shown in \ref{fig:linear-C53T-example}.

Can decks from other SET variants be used to create more variants of this form? Consider the game ProSet/Socks. The group is commutative, and every element is of order 2, which means $ba^{-1}b$ is always equal to $a$. This deck is completely unusable! What about the EvenQuads deck? It can be viewed as $\Z_4^3$, so while some elements have order 2, we can certainly find sets in it. However, there is a different problem with the deck. To play with it, we need to actually assign values to colors and shapes. This means that even if we decide to use this group, we should make different cards! For example, we can style the cards with three squares similar to the pentagons in C53T.

\begin{figure}[ht!]
    \centering
    \begin{tabular}{cc}
        \includegraphics[scale=0.5]{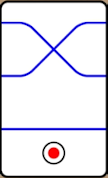} &
        \includegraphics[scale=0.5]{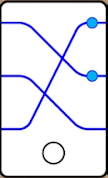} \\
        $a$ & $b$
    \end{tabular}
    \caption{Two cards, $a$ and $b$, showing elements of $C_2\wr S_3$. There is a unique third card that completes an arithmetic set with these two.}
    \label{fig:two-cards-with-beads}
\end{figure}

In Figure~\ref{fig:three-beaded-wires}, we discuss the wreath product $C_2\wr S_3$. Consider the following two cards, seen in the Numberphile video on the variations of the game of SET \cite{N, project}, shown in Figure~\ref{fig:two-cards-with-beads}.

\begin{figure}[ht!]
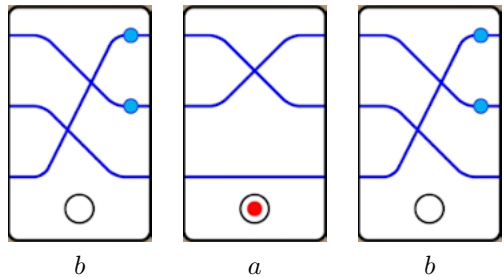

    \centering
    \begin{tabular}{ccc}
        \includegraphics[scale=0.5]{GroupSETbSm.png} &
        \includegraphics[scale=0.5]{GroupSETaSm.png} &
        \includegraphics[scale=0.5]{GroupSETbSm.png} \\
        $b$ & $a$ & $b$
    \end{tabular}
    \caption{The cards $b$, $a$, and $b$. Since $a=a^{-1}$, the product of these is the card $c$, completing the set with $a$ and $b$; this card is equal to $a$, but with additional beads on the first and third wires, in their order on the right side of the card.}
    \label{fig:bab-beads}
\end{figure}

The first card, $a$, is its own inverse, so the third card we are looking for is equal to the product $bab$, which we can visualize as follows. If we ignore the beads, the card that completes the set is $a$. Luckily, if we do not ignore the beads, it is not $a$; we need to add two beads to $a$, to the first and third wires, as discussed in Figure~\ref{fig:bab-beads}.

The products and inverses of permutations are difficult to visualize, especially when the cards are not laid out in a row, so playing this game with the cards in the Numberphile video might be difficult. The good news is that this group can be visualized in many different ways:

\begin{itemize}
    \item The group of symmetries of a cube,
    \item The group of symmetries of an octahedron,
    \item The wreath product of $C_2$ and $S_3$,
    \item The direct product of $C_2$ and $S_4$.
\end{itemize}

\begin{figure}[ht!]
    \centering
\begin{tikzpicture}
  \node[anchor=south west,inner sep=0] (img) at (0,0) {\includegraphics[scale=0.5]{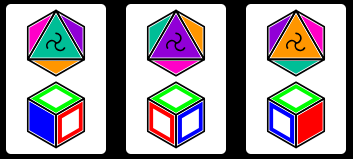}};
  \begin{scope}[shift={(img.south west)},x={(img.south east)},y={(img.north west)}]
    \fill[white,draw=none] (0.32,-0.01) rectangle (0.34,1.01);
    \fill[white,draw=none] (0.66,-0.01) rectangle (0.68,1.01);
  \end{scope}
\end{tikzpicture}
\caption{A set in OCTA Set. There are several ways to interpret each card. We can imagine an octahedron where opposite faces have the same color. To detect whether the octahedron has been reflected, we note the handedness of the spiral shown on it. Similarly, the cube's opposite faces are of the same color. One face in each opposite pair has a white square on it, allowing us to tell if the cube has been reflected. The octahedron can also be seen as the direct product of $S_4$ (the permutation of the four colors in the four triangles) and $C_2$ (the handedness of the spiral). Meanwhile, the cube shows us the wreath product of $C_2$ (whether or not each color is solid or has a white rhombus on it) and $S_3$ (the permutation of the three colors in the rhombi). As a rotation, this set can be seen as rotating about the vertical axis passing through either the green faces of the cube or the top and bottom vertices of the octahedron.}
    \label{fig:octa-set}
\end{figure}

One beautiful deck from the tsetse website \cite{tse} allows one to use any one of the four definitions to play SET with this group. The game is called OCTA Set, as the underlying group is the octahedral group. An example of a set is shown in Figure~\ref{fig:octa-set}, where the top and bottom shapes represent the same element of the group. Moreover, the deck is a torsor: there is no identity card.

\begin{itemize}
    \item The top shape represents an octahedron where the opposite faces are the same color. Thus, the top shape shows a unique rotation of an octahedron. All the swirls are the same in the image, but in other cards, the swirl could be white. Two different swirls mean that the octahedron needs to be reflected to get from one shape to the other. In our example above, no reflection is involved, so all the swirls are black.
    \item We can ignore that the top shape represents an octahedron and only look at the choice of colors. The changes of colors represent a permutation in $S_4$, and a swirl represents an element in $C_2$.
    \item The bottom shape represents a cube where the opposite faces are the same color. One of the faces is solid, and the opposite face is hollow. This way, one can reconstruct the complete coloring of a cube from the top three faces.
    \item We can ignore that the bottom shape represents a cube and only look at the choice of colors. The colors form a permutation, and the beads correspond to changing the hollowness of a color. When viewing the cube this way, each color has its hollowness associated with an element of $C_2$, and the permutation acting on the colors is $S_3$. This exactly corresponds to the wreath product of $C_2$ and $S_3$.
\end{itemize}

Let us prove that this is a set. Consider the bottom cube shape. Comparing the first two cards, the top face does not change. We can see that the symmetry of the cube is the 90-degree clockwise rotation around the line that goes through the centers of the green faces. In such a rotation, the left face on the second card keeps the color from the first card, while the right face takes the color from the left face on the first card and swaps hollowness. We see that the third cube completes the set.

For another proof, let us look at the top shape and discuss what happens with the permutation of colors when changing from the first card to the second. The left color moves to the bottom, the bottom color to the right, the right color to the center, and the center color to the left. Not surprisingly, we got a cyclic permutation of order 4, similar to a 90-degree rotation being of order 4. The same permutation happens when moving from the second card to the third. The swirl stays the same.

It is worth emphasizing that the two shapes on an OCTA Set card are just two different illustrations of the same group element.
It is thus possible to look only at the octahedra (or only at the cubes) to play OCTA Set. In practice, certain sets are easier to spot when looking at the octahedra, while others are easier to see on the cubes.
Furthermore, as the polyhedra are dual to each other, they have the same isometries (rotations and reflections).
This is reflected on the cards, so any set found has the same transformation to its polyhedra.
For example, in Figure~\ref{fig:octa-set}, both polyhedra are rotated about the vertical axis.

As we mentioned, when you play this game and pick two cards and then calculate what card completes the set, you might discover that it is one of the cards you picked. The probability that two cards in a specific order cannot be completed into a set is the same as the probability of picking a random element in our group and discovering that it has order 2. The symmetric group $S_4$ has 9 elements of order 2. Thus, the direct product with $S_2$ has 19 elements of order 2, giving a probability of 19/48. For completeness, this group also has 8 elements of order 3, 12 elements of order 4, and 8 elements of order 6, not to mention the identity of order 1.

If $ba^{-1}$ is an element of order three, then the cards $a$, $b$, and the card $c$ that completes the set form a set when they are taken in any order. This can be seen in the game of SET, the group element $ba^{-1}$ always has order 3.

Notably, there was nothing special or extraordinary about the group discussed above. It has a pretty visualization as a cube or octahedron, but is not otherwise particularly interesting. The reason why this group allowed for these two platonic solids to be used to visualize it is because the cube is dual to the octahedron. But we could have similarly used any group to play SET! One such example might consider using the group of rotations of the other pair of dual Platonic solids, the icosahedron and the dodecahedron. This group is actually equivalent to $A_5$, also known as the alternating group of order 5, which consists of all even permutations of five elements. The tsetse website \cite{tse} contains an implementation of such a game called A5SET (pronounced ``asset''). The design of the site, games, and cards was done by Andy Tockman and Della Hendrickson.

\section{What Groups make Good Set Games}
\label{sec:good}

There are several considerations for making a good game out of a group.

\textbf{The size of the deck.} The size of the deck is the order of the group, if the identity is included. We want the deck to be playable. That means it should not be too small or too big, as the first will run out too quickly and the second is too hard to shuffle and will play too slowly. The SET deck has 81 cards.

\textbf{The number of cards on the table.} The rules of the game can customize the number of cards that are usually on the table. The number of visible cards can be increased if there is no set. Having too many visible cards makes the gameplay too overwhelming, but having too few cards will frequently result in having to add more. The game of SET starts with 12 cards. The probability of a set to exist among 12 cards is approximately 96.77\%, relatively high. In the game of ProSet, seven cards are guaranteed to have a set among them. In EvenQuads, 10 cards are needed to guarantee a set~\cite{CragerEtAl}. In a good game, we do not want to start with much more than a dozen cards, but we want to have a reasonably high probability of finding a set.

\textbf{Visual presentation.} The elements of the group need to have a nice visual representation. If there are several group descriptions that are isomorphic to the chosen group, it can be useful to depict each element of several groups on a card, as in OCTA Set. Aside from aesthetic reasons, it can be the case that finding some sets is easier by examining one group presentation, while others might be easier by examining a different one. Additionally, when there is an elegant property uniting sets, a visual presentation can make this easier to see, for example, when playing C53T, sets have an axis of symmetry across each of the three pentagons, but this would not be visible if the cards had triplets of numbers from 0 to 4!

\textbf{Torsor property.} It is good when the cards are equivalent in a sense, as for example, in SET or EvenQuads. This happens when we have a torsor. The Socks deck is not based on a torsor, and the identity card is removed. In this game, when there is a small number of socks on cards, finding sets is easier than otherwise.

\textbf{The number of cards in a set.} In games based on arithmetic progressions in groups, sets will always have three cards. However, in a game such as C53T, where sets have 5 cards, finding a set is incredibly difficult!

\textbf{The number of cards added when no sets exist.} If there are no sets in the visible cards, some cards should be added to the table. In SET, three cards are added at a time when the visible cards contain no sets. The number of added cards does not have to be three, but if only one card is added at a time, then any new set is guaranteed to contain the newly added card, making sets much easier to find! Adding new cards one at a time is thus potentially an undesirable feature of a set game.

Armed with this knowledge, we have assembled all the ingredients needed to create even more SET variants and introduce an ever-growing number of people to this amazing world!

\end{document}